\documentclass[11pt]{article}
\usepackage[ansinew]{inputenc}
\usepackage{array}
\usepackage{color}
\usepackage{amsmath,amsthm,amscd}
\usepackage{amsxtra}
\usepackage{amstext}
\usepackage{amssymb}
\usepackage{hyperref}
\usepackage{latexsym}
\usepackage{amsfonts,cite}
\usepackage{graphics,epstopdf}
\usepackage{epsfig}
\usepackage{amsfonts}
\usepackage{graphicx}
\usepackage{tcolorbox}
\usepackage{array,tabularx,colortbl}
\usepackage[framemethod=TikZ]{mdframed}

\providecommand{\U}[1]{\protect\rule{.1in}{.1in}}

\begin{document}
	
	\sloppy
	\newtheorem{thm}{Theorem}
	\newtheorem{cor}{Corollary}
	\newtheorem{lem}{Lemma}
	\newtheorem{prop}{Proposition}
	\newtheorem{eg}{Example}
	\newtheorem{defn}{Definition}
	\newtheorem{rem}{Remark}
	\newtheorem{note}{Note}
	\numberwithin{equation}{section}
	
	\thispagestyle{empty}
	\parindent=0mm
	
	\begin{center}
		{\large \textbf{Umbral insights into a hybrid family of hypergeometric\\
				\vspace{0.15cm} and Mittag-Leffler functions}}\\ 					
		
		\vspace{0.20cm}
		
		{\small\bf Subuhi Khan$^{1}$, Ujair Ahmad$^{2}$, Mehnaz Haneef$^{3}$}\\
		\vspace{0.15cm}
		
		Department of Mathematics,\\ Aligarh Muslim University, Aligarh-202001, U.P., India \\
		\footnote{$^{1}$E-mail:~subuhi2006@gmail.com (Subuhi Khan) (Corresponding author)}
		\footnote{$^{2}$E-mail:~gi8986@myamu.ac.in (Ujair Ahmad)}
		\footnote{$^{3}$E-mail:~mehnaz272@gmail.com (Mehnaz Haneef)}
		
	\end{center}
	
	\begin{abstract}
		\noindent
		The umbral approach provides methods for comprehending and redefining special functions. This approach is employed efficiently in order to uncover intricacies and introduce new families of special functions. In this article, the umbral perspective is adopted to introduce a hybrid family of hypergeometric and Mittag-Leffler functions. The umbral-operational procedures are used to derive the generating functions, explicit representations, differential recurrence formulae, and specific integral formulae. Further, the Laplace and Sumudu transforms for the hypergeometric-Mittag-Leffler functions are established. The graphical representation and pattern for distribution of zeros for suitable values of parameters are also presented.	
	\end{abstract}
	\parindent=0mm
	\vspace{.25cm}
	
	\noindent
	\textbf{Key Words:}~~Umbral methods; Hypergeometric functions; Mittag-Leffler functions; Hypergeometric-Mittag-Leffler functions.
	
	\vspace{0.25cm}
	\noindent
	\textbf{2020 Mathematics Subject Classification:}~~05A40; 33C05; 33C20; 33E12; 44A10; 44A99.
	
	\section{Preliminaries}
	Hypergeometric functions \cite{AND} are generalizations of the higher-order transcendental functions, captivating with the intriguing properties that make them invaluable. The theoretical framework of hypergeometric functions has been recently reformulated \cite{DMS} based on the indicial umbral calculus (IUC), developed in the last two decades by Dattoli and coworkers, see for example \cite{DSLM,DGH} (for a recent updated review see \cite{LD}). The IUC formalism deepens its roots in the classical Roman-Rota umbral calculus \cite{RT}, which fixes the abstract notation and the operational rules allowing the handling of formal series.\\
	
	The IUC is an evolution of the monomiality principle \cite{ST}, which in turn share analogies with the operational calculus. The IUC allows the possibility of formal handling of higher transcendental functions as ordinary transcendent or rational functions. The formal tools to accomplish such a realization are created through the ``image" and the ``vacuum" functions. Within this framework, the image of the Bessel functions is a Gaussian or a Lorentzian, depending on the adopted vacuum function.\\

The function
\begin{equation*}
	\zeta_{u}:=\dfrac{1}{\Gamma(u+1)},
\end{equation*}
is called the vacuum function or umbral vacuum. The umbral operator or shift operator $\hat{c}$ is defined by:
\begin{equation*}
	\hat{c}=e^{\partial u},
\end{equation*}
with {\em u} being the variable on which the operator acts.\\
We recall that $\hat{c}^{\mu}$ acts on the corresponding vacuum $\zeta_{0}$ as follows:
\begin{equation}\label{mheq10}
	\hat{c}^{\mu}\zeta_{0}=\dfrac{1}{\Gamma(\mu+u+1)}\bigg|_{u=0}=\dfrac{1}{\Gamma(\mu+1)},\quad\forall\mu\in\mathbb{R}.
\end{equation}
\noindent  The hypergeometric functions are encased within umbral formalism for the first time by Dattoli {\em et al.} \cite{DMS} by introducing their umbral form as:
\begin{equation}\label{mheq5}
	{_{2}F_{1}(\mathfrak{a_{1},a_{2};b_{1};}u)}=e^{u\;\hat{_{2}\chi_{1}}}\Phi_{0}.
\end{equation}
Here, the vacuum function
\begin{equation*}
	\Phi_{r}=\dfrac{(\mathfrak{a_{1}})_{r}(\mathfrak{a_{2}})_{r}}{(\mathfrak{b_{1}})_{r}},
\end{equation*}
with the corresponding umbral operator or vacuum shift operator
\begin{equation*}
	\hat{_{2}\chi_{1}}=e^{\partial u_{1}}e^{\partial u_{2}}e^{\partial u_{3}},
\end{equation*}
acts on vacuum $\Phi_{0}$ such that
\begin{equation}\label{mheq6}
	_{2}\hat{\chi}^{r}_{1}\Phi_{0}:=\Phi_{r}=\dfrac{(\mathfrak{a_{1}})_{r}(\mathfrak{a_{2}})_{r}}{(\mathfrak{b_{1}})_{r}}.
\end{equation}
In view of equations \eqref{mheq5} and \eqref{mheq6}, the series representation for hypergeometric functions ${_{2}F_{1}(\mathfrak{a_{1},a_{2};b_{1};}u)}$ is obtained.

\begin{note}
	The vacuum function can be specified by the following function:
	\begin{equation*}
		\Phi(\mathfrak{a_{1}},u_{1},\mathfrak{a_{2}},u_{2},\mathfrak{b_{1}},u_{3})=\dfrac{(\mathfrak{a_{1}})_{u_{1}}(\mathfrak{a_{2}})_{u_{2}}}{(\mathfrak{b_{1}})_{u_{3}}}.
	\end{equation*}
	with $u_{1}$, $u_{2}$, $u_{3}$ as the domain's variables of the function on which the operator acts.\\
	Further
	\begin{equation*}
		e^{d\partial u}(a)_{u}=\dfrac{\Gamma(u+d+a)}{\Gamma(a)}.
	\end{equation*}
	In conclusion, we can write
	\begin{equation*}
		e^{r\;\partial u_{1}}e^{r\;\partial u_{2}}e^{r\;\partial u_{3}}	\Phi(\mathfrak{a_{1}},u_{1},\mathfrak{a_{2}},u_{2},\mathfrak{b_{1}},u_{3})|_{u_{1}=u_{2}=u_{3}=0}= {_{2}\hat{\chi}^{r}_{1}}\Phi_{0}:=\dfrac{(\mathfrak{a_{1}})_{r}(\mathfrak{a_{2}})_{r}}{(\mathfrak{b_{1}})_{r}}.
	\end{equation*}
\end{note}
Further, we recall that the generalized hypergeometric functions are defined by the following seires \cite{AND}:
\begin{equation}\label{mheq4}
	_{p}F_{q}(\mathfrak{a_{1},a_{2},\cdots,a}_{p};\mathfrak{b_{1},b_{2},\cdots,b}_{q};u)=\sum_{r=0}^{\infty}\dfrac{\mathfrak{(a_{1})}_{r}\mathfrak{(a_{2})}_{r}\cdots\mathfrak{(a_{p})}_{r}}{\mathfrak{(b_{1})}_{r}\mathfrak{(b_{2})}_{r}\cdots\mathfrak{(b_{q})}_{r}}\dfrac{u^{r}}{r!}.
\end{equation}
The umbral form of the generalized hypergeometric functions is introduced as \cite{DMS}:
\begin{equation}\label{mheqi}
	_{p}F_{q}(\mathfrak{a_{1},a_{2},\cdots,a}_{p};\mathfrak{b_{1},b_{2},\cdots,b}_{q};u)=e^{u\;\hat{_{p}\chi_{q}}}\Phi_{0},
\end{equation}
where, $_{p}\hat{\chi}_{q}$ is the umbral operator which acts on the umbral vacuum $\Phi_{0}$ such that 
\begin{equation}\label{mheqii}
	_{p}\hat{\chi}^{r}_{q}\Phi_{0}=\sum_{r=0}^{\infty}\dfrac{\mathfrak{(a_{1})}_{r}\mathfrak{(a_{2})}_{r}\cdots\mathfrak{(a_{p})}_{r}}{\mathfrak{(b_{1})}_{r}\mathfrak{(b_{2})}_{r}\cdots\mathfrak{(b_{q})}_{r}}.
\end{equation}
\noindent We recall that the classical Mittag-Leffler functions of one parameter were initially studied in \cite{M} to address divergent series summation methods. These functions play pivotal role in solving differential and integral equations of fractional order \cite{KST}. These equations emerge in various applied scientific disciplines such as physics, chemistry, engineering, and biology. The Mittag-Leffler functions hold a fundamental significance in fractional calculus, akin to the exponential function in traditional calculus, driving its extensive study.\\

Wiman \cite{W} introduced the 2-parameter Mittag-Leffler functions in the following form: 
\begin{equation}\label{mheq51}
E_{\alpha,\beta}(u)=\sum_{m=0}^{\infty}\dfrac{u^{m}}{\Gamma(\alpha m+\beta)},\quad\quad u\in\mathbb{R};\quad\quad\alpha,\beta\in\mathbb{R^{+}}.
\end{equation}
The Mittag-Leffler functions play a significant role in the umbral calculus. Dattoli {\em et al.} \cite{DGH} considered the umbral restyling of the Mittag-Leffler functions using operator \eqref{mheq10} in the following form:
\begin{equation}\label{mheq9}
E_{\alpha,\beta}(u)=\dfrac{\hat{c}^{(\beta-1)}}{1-\hat{c}^\alpha u}\zeta_{0}.
\end{equation}
Another umbral operator $\hat{d}_{(\alpha,\beta)}$, which acts on the vacuum function 
\begin{equation}\label{mheqviii}
\Psi_{u}:=\dfrac{\Gamma(u+1)}{\Gamma(\alpha u+\beta)},
\end{equation}
is taken as \cite{DGH}:
\begin{equation}\label{mheqx}
\hat{d}^{k}_{(\alpha,\beta)}\Psi_{u}=\dfrac{\Gamma(k+u+1)}{\Gamma(\alpha(k+u)+\beta)},\quad\quad k\in\mathbb{R};\quad\quad k+\mathbf{Re}(u)\ge-1.
\end{equation}
From equation \eqref{mheqx}, it follows that
\begin{equation}\label{mheq8}
\hat{d}^{k}_{(\alpha,\beta)}\Psi_{0}:=\hat{d}^{k}_{(\alpha,\beta)}\Psi_{u}\bigg|_{u=0}=\dfrac{\Gamma(k+1)}{\Gamma(\alpha k+\beta)},\quad\quad k\in\mathbb{R};\quad\quad k\ge-1.
\end{equation}
Consequently, in view of equations \eqref{mheq51} and \eqref{mheq8}, an alternate umbral form of the $2$-parameter Mittag-Leffler functions is given as \cite{DGH}:
\begin{equation}\label{mheq7}
E_{\alpha,\beta}(u)=e^{u\;\hat{d}_{(\alpha,\beta)}}\Psi_{0}.
\end{equation}	
\noindent 
Inspired by the applications of umbral methods, this article aims to introduce a distinctive hybrid family of special functions termed as hypergeometric-Mittag-Leffler functions.\\

\noindent The layout of this article unfolds as follows: In Section $2$, the hypergeometric-Mittag-Leffler functions are introduced using umbral techniques. The series expansion, generating equation and certain differential relations are derived. In Section $3$, integral representations for the hypergeometric-Mittag-Leffler functions are derived. In Section $4$, the Laplace and Sumudu transforms for the hypergeometric-Mittag-Leffler functions along with graphical illustrations of functions for certain specific values of parameters are explored. The relations of the hypergeometric-Mittag-Leffler functions with other functions are considered in Appendix.\\
\section{Hypergeometric-Mittag-Leffler functions}	
To introduce the hypergeometric-Mittag-Leffler functions (HMLF) denoted by ${_{_{2}F_{1}}{E_{(\alpha,\beta)}}}(\mathfrak{a_{1},a_{2};b_{1};}u)$, we make use of the umbral operator of hypergeometric functions $_{2}F_{1}(\mathfrak{a_{1};a_{2};b_{1};}u)$ in the Mitag-Leffler functions ${E_{\alpha,\beta}}(u)$. Replacing the variable $u$ in ${E_{\alpha,\beta}}(u)$ by the umbral operator $_{2}\hat{\chi}_{1}$ times $u$, we have
\begin{equation}\label{mheq11}
{_{_{2}F_{1}}{E_{(\alpha,\beta)}}}(\mathfrak{a_{1},a_{2};b_{1};}u)=E_{\alpha,\beta}(u\;\hat{_{2}\chi_{1}})\Phi_{0}.
\end{equation} 
Using umbral form \eqref{mheq9} of Mittag-Leffler functions in r.h.s. of equation \eqref{mheq11}, the following umbral form of the HMLF ${_{_{2}F_{1}}{E_{(\alpha,\beta)}}}(\mathfrak{a_{1},a_{2};b_{1};}u)$ is obtained in terms of the operator $\hat{c}$:
\begin{equation}\label{mheq13}
{_{_{2}F_{1}}{E_{(\alpha,\beta)}}}(\mathfrak{a_{1},a_{2};b_{1};}u)=\dfrac{\hat{c}^{(\beta-1)}}{1-\hat{c}^{\alpha}\;u\;{_{2}\hat{\chi}_{1}}}\Phi_{0}\zeta_{0}.
\end{equation}
Further, using umbral form \eqref{mheqi} of the generalized hypergeometric functions and umbral operator \eqref{mheqii}, the umbral form of the generalized HMLF ${_{_{p}F_{q}}{E_{(\alpha,\beta)}}}(\mathfrak{a_{1},a_{2},\cdots,a}_{p};\mathfrak{b_{1},b_{2},\cdots,b}_{q};u)$ in terms of the operator $\hat{c}$ is given as:
\begin{equation}\label{mheqiv}
{_{_{p}F_{q}}{E_{(\alpha,\beta)}}}(\mathfrak{a_{1},a_{2},\cdots,a}_{p}\mathfrak{;b_{1},b_{2},\cdots,b}_{q};u)=\dfrac{\hat{c}^{(\beta-1)}}{1-\hat{c}^{\alpha}\;u\;{_{p}\hat{\chi}_{q}}}\Phi_{0}\zeta_{0}.
\end{equation}
Again, by using the alternate umbral form \eqref{mheq7} of Mittag-Leffler functions in equation \eqref{mheq11}, the following alternate umbral form of HMLF ${_{_{2}F_{1}}{E_{(\alpha,\beta)}}}(\mathfrak{a_{1},a_{2};b_{1};}u)$ is obtained:
\begin{equation}\label{mheq12}
{_{_{2}F_{1}}{E_{(\alpha,\beta)}}}(\mathfrak{a_{1},a_{2};b_{1};}u)=e^{u\;{_{2}\hat{\chi}_{1}}\;\hat{d}_{(\alpha,\beta)}}\Psi_{0}\Phi_{0}.
\end{equation}
In view of equation \eqref{mheqi} and operator \eqref{mheqii}, the umbral form of the generalized HMLF ${_{_{p}F_{q}}{E_{(\alpha,\beta)}}}(\mathfrak{a_{1},a_{2},\cdots,a}_{p}\mathfrak{;b_{1},b_{2},\cdots,b}_{q};u)$ is given as:
\begin{equation}\label{mheqv}
{_{_{p}F_{q}}{E_{(\alpha,\beta)}}}(\mathfrak{a_{1},a_{2},\cdots,a}_{p}\mathfrak{;b_{1},b_{2},\cdots,b}_{q};u)=e^{u\;{_{p}\hat{\chi}_{q}}\;\hat{d}_{(\alpha,\beta)}}\Psi_{0}\Phi_{0}.
\end{equation}
Next, we establish the series representation of the HMLF ${_{_{2}F_{1}}{E_{(\alpha,\beta)}}}(\mathfrak{a_{1},a_{2};b_{1};}u)$ by proving the following result:
\begin{thm}
For the hypergeometric-Mittag-Leffler functions ${_{_{2}F_{1}}{E_{(\alpha,\beta)}}}(\mathfrak{a_{1},a_{2};b_{1};}u)$, the following series expansion holds true:
\begin{equation}\label{mheq14}
	{_{_{2}F_{1}}{E_{(\alpha,\beta)}}}(\mathfrak{a_{1},a_{2};b_{1};}u)=\sum_{r=0}^{\infty}\dfrac{(\mathfrak{a}_{1})_{r}(\mathfrak{a}_{2})_{r}}{(\mathfrak{b}_{1})_{r}\Gamma(\alpha r+\beta)}u^{r},\quad\quad\alpha,\beta\in\mathbb{R^{+}},\quad\quad u\in\mathbb{R}.
\end{equation}
\begin{proof} Expanding the denominator in equation \eqref{mheq13} and simplifying, it follows that
	\begin{equation}\label{mheqvii}
		{_{_{2}F_{1}}{E_{(\alpha,\beta)}}}(\mathfrak{a_{1},a_{2};b_{1};}u)	=\sum_{r=0}^{\infty}\hat{c}^{(\alpha r+\beta-1)}\;u^{r}\;{_{2}\hat{\chi}^{r}_{1}}\Phi_{0}\zeta_{0},
	\end{equation}
	which on using operator actions \eqref{mheq10} and \eqref{mheq6} yields assertion \eqref{mheq14}.	\\
	
	Alternately, expanding the exponetial function in the r.h.s. of equation \eqref{mheq12}, we find
	\begin{equation}\label{mheq15}
		{_{_{2}F_{1}}{E_{(\alpha,\beta)}}}(\mathfrak{a_{1},a_{2};b_{1};}u)=\sum_{r=0}^{\infty}\dfrac{u^{r}\hat{_{2}\chi^{r}_{1}}\;\hat{d}^{r}_{(\alpha,\beta)}}{r!}\Psi_{0}\Phi_{0}.
	\end{equation}
	Use of operator actions \eqref{mheq6} and \eqref{mheq8} in equation \eqref{mheq15} also leads to assertion \eqref{mheq14}, which confirms the fact that umbral forms \eqref{mheq13} and \eqref{mheq12} of HMLF ${_{_{2}F_{1}}{E_{(\alpha,\beta)}}}(\mathfrak{a_{1},a_{2};b_{1};}u)$ are equivalent.
	
\end{proof}
\end{thm}
Considering umbral form \eqref{mheqv} of the generalized HMLF ${_{_{p}F_{q}}{E_{(\alpha,\beta)}}}(\mathfrak{a_{1},a_{2},\cdots,a}_{p}\mathfrak{;b_{1},b_{2},\cdots,b}_{q};u)$ and following the lines proof of Theorem \thmnumber{1}, we obtain the following result:
\begin{thm}
For the generalized hypergeometric-Mittag-Leffler functions ${_{_{p}F_{q}}{E_{(\alpha,\beta)}}}(\mathfrak{a_{1},a_{2},\cdots,a}_{p}\mathfrak{;b_{1},b_{2},\cdots,b}_{q};u)$, the following series expansion holds true:
\begin{equation}\label{mheq86}
	{_{_{p}F_{q}}{E_{(\alpha,\beta)}}}(\mathfrak{a_{1},a_{2},\cdots,a}_{p}\mathfrak{;b_{1},b_{2},\cdots,b}_{q};u)=\sum_{r=0}^{\infty}\dfrac{\mathfrak{(a_{1})}_{r}\mathfrak{(a_{2})}_{r}\cdots\mathfrak{(a_{p})}_{r}}{\mathfrak{(b_{1})}_{r}\mathfrak{(b_{2})}_{r}\cdots\mathfrak{(b_{q})}_{r}}\dfrac{u^{r}}{\Gamma(\alpha r+\beta)}.
\end{equation}
If $p\leq q$, the series converges for all finite $u$ and if $p=q+1$, the series converges for $|u|< 1$ and diverges for $|u|>1$. If $p>q+1$, the series diverges for $u\not=0$.
\end{thm}
\noindent	Next, we establish an identity for the HMLF ${_{_{2}F_{1}}{E_{(\alpha,\beta)}}}(\mathfrak{a_{1},a_{2};b_{1};}u)$ by proving the following result:
\begin{thm} For $\alpha,\beta\in\mathbb{R^{+}}$ and $n\in\mathbb{N}$, the following identity involving hypergeometric-Mittag-Leffler functions ${_{_{2}F_{1}}{E_{(\alpha,\beta)}}}(\mathfrak{a_{1},a_{2};b_{1};}u)$ holds:
\begin{align}\label{mheq16}
	\notag
	{_{_{2}F_{1}}{E_{(\alpha,\beta-n\alpha)}}}(\mathfrak{a_{1},a_{2};b_{1};}u)=&\dfrac{u^{n}(\mathfrak{a_{1}})_{n}(\mathfrak{a_{2}})_{n}}{(\mathfrak{b_{1}})_{n}}{_{_{2}F_{1}}{E_{(\alpha,\beta)}}}(\mathfrak{a_{1}}+n,\mathfrak{a_{2}}+n;\mathfrak{b_{1}}+n;u)\\
	&+\sum_{r=1}^{n}\dfrac{(\mathfrak{a_{1}})_{n-r}(\mathfrak{a_{2}})_{n-r}}{(\mathfrak{b_{1}})_{n-r}}\dfrac{u^{(n-r)}}{\Gamma(\beta-r\alpha)};\quad \beta>r\alpha, \;\forall r\leq n\in\mathbb{N}.
\end{align}
\begin{proof}
	Replacing $\beta$ by $\beta-\alpha$ in equation \eqref{mheq14}, it follows that: 
	\begin{equation*}
		{_{_{2}F_{1}}{E_{(\alpha,\beta-\alpha)}}}(\mathfrak{a_{1},a_{2};b_{1};}u)=\dfrac{1}{\Gamma(\beta-\alpha)}+\sum_{r=1}^{\infty}\dfrac{(\mathfrak{a_{1}})_{r}(\mathfrak{a_{2}})_{r}}{(\mathfrak{b_{1}})_{r}\Gamma(\alpha (r-1)+\beta)}u^{r},
	\end{equation*}
	which on replacing $r$ by $s+1$ in the r.h.s. and rearranging the series takes the form
	\begin{equation}\label{mheq17}
		{_{_{2}F_{1}}{E_{(\alpha,\beta-\alpha)}}}(\mathfrak{a_{1},a_{2};b_{1};}u)=\dfrac{1}{\Gamma(\beta-\alpha)}+\dfrac{u\mathfrak{a_{1}a_{2}}}{\mathfrak{b_{1}}}\sum_{s=0}^{\infty}\dfrac{(\mathfrak{a_{1}}+1)_{s}(\mathfrak{a_{2}}+1)_{s}}{(\mathfrak{b_{1}}+1)_{s}\Gamma(\alpha s+\beta)}u^{s}.
	\end{equation}
	Using series definition \eqref{mheq14} in the r.h.s. of equation \eqref{mheq17}, we find
	\begin{equation}\label{mheq83}
		{_{_{2}F_{1}}{E_{(\alpha,\beta-\alpha)}}}(\mathfrak{a_{1},a_{2};b_{1};}u)=\dfrac{u\mathfrak{a_{1}a_{2}}}{\mathfrak{b_{1}}}{_{_{2}F_{1}}{E_{(\alpha,\beta)}}}(\mathfrak{a_{1}}+1,\mathfrak{a_{2}}+1;\mathfrak{b_{1}}+1;u)+\dfrac{1}{\Gamma(\beta-\alpha)}.
	\end{equation}
	Next, replacing $\beta$ by $\beta-2\alpha$ in equation \eqref{mheq14}, we find
	\begin{align*}
		{_{_{2}F_{1}}{E_{(\alpha,\beta-2\alpha)}}}(\mathfrak{a_{1},a_{2};b_{1};}u)=\dfrac{1}{\Gamma(\beta-2\alpha)}+\dfrac{u\mathfrak{a_{1}a_{2}}}{\mathfrak{b_{1}}\Gamma(\beta-\alpha)}+\sum_{r=2}^{\infty}\dfrac{(\mathfrak{a_{1}})_{r}(\mathfrak{a_{2}})_{r}}{(\mathfrak{b_{1}})_{r}\Gamma(\alpha (r-2)+\beta)}u^{r}.
	\end{align*}
	Replacing $r$ by $s+2$ in the r.h.s. and following the same procedure as above, we find
	\begin{equation}\label{mheq84}
		{_{_{2}F_{1}}{E_{(\alpha,\beta-2\alpha)}}}(\mathfrak{a_{1},a_{2};b_{1};}u)	=\dfrac{u^{2}\mathfrak{(a_{1})_{2}(a_{2})_{2}}}{\mathfrak{(b_{1})_{2}}}{_{_{2}F_{1}}{E_{(\alpha,\beta)}}}(\mathfrak{a_{1}}+2,\mathfrak{a_{2}}+2;\mathfrak{b_{1}}+2;u)+\dfrac{1}{\Gamma(\beta-2\alpha)}+\dfrac{u\mathfrak{a_{1}a_{2}}}{\mathfrak{b_{1}}\Gamma(\beta-\alpha)}.
	\end{equation}
	Continuing the process of iteration, assertion \eqref{mheq16} is proved.
\end{proof}
\end{thm}
\noindent Further, we derive certain differential relations for the HMLF ${_{_{2}F_{1}}{E_{(\alpha,\beta)}}}(\mathfrak{a_{1},a_{2};b_{1};}u)$ in terms of the generalized HMLF ${_{_{p}F_{q}}{E_{(\alpha,\beta)}}}(\mathfrak{a_{1},a_{2},\cdots,a}_{p}\mathfrak{;b_{1},b_{2},\cdots,b}_{q};u)$ in the form of following result:
\begin{thm}
For the hypergeometric-Mittag-Leffler functions ${_{_{2}F_{1}}{E_{(\alpha,\beta)}}}(\mathfrak{a_{1},a_{2};b_{1};}u)$, the following differential relations hold true:
\begin{equation}\label{mheq23}
D_{u}^{m}[{_{_{2}F_{1}}{E_{(\alpha,\beta)}}}(\mathfrak{a_{1},a_{2};b_{1};}u)]=\dfrac{m!(\mathfrak{a_{1}})_{m}(\mathfrak{a_{2}})_{m}}{(\mathfrak{b_{1}})_{m}}{_{_{3}F_{2}}{E_{(\alpha,m\alpha+\beta)}}}(\mathfrak{a_{1}}+m,\mathfrak{a_{2}}+m,m+1;\mathfrak{b_{1}}+m,1;u).
\end{equation}
\begin{equation}\label{mheq24}
D_{\hat{d}_{(\alpha,\beta)}}^{m}[e^{u\;{_{2}\hat{\chi}_{1}}\;\hat{d}_{(\alpha,\beta)}}\Psi_{0}\Phi_{0}]=\dfrac{u^{m}(\mathfrak{a_{1}})_{m}(\mathfrak{a_{2}})_{m}}{(\mathfrak{b_{1}})_{m}}{_{_{2}F_{1}}{E_{(\alpha,\beta)}}}(\mathfrak{a_{1}}+m,\mathfrak{a_{2}}+m;\mathfrak{b_{1}}+m;u).\quad\quad\quad\quad
\end{equation}
\begin{equation}\label{mheq25}
\prod_{n=1}^{m}D_{_{n+1}\hat{\chi}_{n}}[e^{u\;{_{2}\hat{\chi}_{1}}\;\hat{d}_{(\alpha,\beta)}}\Psi_{0}\Phi_{0}]=u^{m}{_{_{m+2}F_{m+1}}{E_{(\alpha,m\alpha+\beta)}}}(\mathfrak{a_{1},a_{2}},2,2,\cdots,2;\mathfrak{b_{1}},1,1,\cdots,1;u).\quad\quad
\end{equation}
\begin{proof}
Differentiating equation \eqref{mheq12} w.r.t. $u$, we find
\begin{equation*}
	D_{u}[{_{_{2}F_{1}}{E_{(\alpha,\beta)}}}(\mathfrak{a_{1},a_{2};b_{1};}u)]={_{2}\hat{\chi}_{1}}\;\hat{d}_{(\alpha,\beta)}e^{u\;{_{2}\hat{\chi}_{1}}\;\hat{d}_{(\alpha,\beta)}}\Psi_{0}\Phi_{0},
\end{equation*}
Expanding the exponential and rearranging the terms, it follows that
\begin{equation*}
	D_{u}[{_{_{2}F_{1}}{E_{(\alpha,\beta)}}}(\mathfrak{a_{1},a_{2};b_{1};}u)]=\sum_{r=0}^{\infty}\dfrac{u^{r}\;{_{2}\hat{\chi}^{(r+1)}_{1}}\;\hat{d}^{(r+1)}_{(\alpha,\beta)}}{r!}\Psi_{0}\Phi_{0}.
\end{equation*}
Applying operator equations \eqref{mheq6} and \eqref{mheq8} and simplifying the resultant equation by making use of identity
\begin{equation}\label{mheqxi}
	(d)_{n+m}=(d)_{m}(d+m)_{n}\;,
\end{equation}
we have
\begin{equation*}
	D_{u}[{_{_{2}F_{1}}{E_{(\alpha,\beta)}}}(\mathfrak{a_{1},a_{2};b_{1};}u)]=\sum_{r=0}^{\infty}\dfrac{\mathfrak{a_{1}}(\mathfrak{a_{1}}+1)_{r}\mathfrak{a_{2}}(\mathfrak{a_{2}}+1)_{r}(r+1)\Gamma(r+1)}{\mathfrak{b_{1}}(\mathfrak{b_{1}}+1)_{r}r!\Gamma(\alpha r+\alpha+\beta)}u^{r}.
\end{equation*}
Expressing $r+1$ as $\dfrac{(2)_{r}}{(1)_{r}}$ in above equation and using definition \eqref{mheq86} for $p=3$ and $q=2$, we find
\begin{equation*}
	D_{u}[{_{_{2}F_{1}}{E_{(\alpha,\beta)}}}(\mathfrak{a_{1},a_{2};b_{1};}u)]=\dfrac{\mathfrak{a_{1}a_{2}}}{\mathfrak{b_{1}}}{_{_{3}F_{2}}{E_{(\alpha,\alpha+\beta)}}}(\mathfrak{a_{1}}+1,\mathfrak{a_{2}}+1,2;\mathfrak{b_{1}}+1,1;u),
\end{equation*}
which shows that result \eqref{mheq23} holds true for $m=1$.\\
Assuming that result \eqref{mheq23} holds for $m=k$, so that in view of equations \eqref{mheqv} and \eqref{mheq86}, we have
\begin{equation}\label{mheq61}
	D_{u}^{k}[{_{_{2}F_{1}}{E_{(\alpha,\beta)}}}(\mathfrak{a_{1},a_{2};b_{1};}u)]=\dfrac{k!(a)_{k}(b)_{k}}{(c)_{k}}e^{u\;{_{3}\hat{\chi}_{2}}\;\hat{d}_{(\alpha,k\alpha+\beta)}}\Psi_{0}\Phi_{0},
\end{equation}
where
\begin{equation*}
	{_{3}\hat{\chi}^{r}_{2}}\Phi_{0}=\dfrac{\mathfrak{(a_{1}}+k)_{r}\mathfrak{(a_{2}}+k)_{r}(k+1)_{r}}{\mathfrak{(b_{1}}+k)_{r}(1)_{r}}.
\end{equation*}
Differentiating equation \eqref{mheq61} w.r.t. $u$, it gives
\begin{equation}\label{mheq95}
	D_{u}^{k+1}[{_{_{2}F_{1}}{E_{(\alpha,\beta)}}}(\mathfrak{a_{1},a_{2};b_{1};}u)]=\dfrac{k!(a)_{k}(b)_{k}}{(c)_{k}}\sum_{r=0}^{\infty}\dfrac{u^{r}\;{_{3}\hat{\chi}^{r+1}_{2}}\;\hat{d}_{(\alpha,k\alpha+\beta)}^{r+1}}{r!}\Psi_{0}\Phi_{0}.
\end{equation}
Applying operator equations \eqref{mheqii} and \eqref{mheq8} in the r.h.s. of equation \eqref{mheq95} and using definition \eqref{mheq86}, it follows that equation \eqref{mheq23} holds true for $m=k+1$. Hence by induction, assertion \eqref{mheq23} is proved.\\

\noindent Proceeding on the same lines as above, assertions \eqref{mheq24} and \eqref{mheq25} can be proved.
\end{proof}
\end{thm}
\noindent In order to further bolster the legitimacy of the umbral formalism, we delve into certain integral representations in the next section.
\section{Integral representations}
In this section, certain integral representation for the HMLF ${_{_{2}F_{1}}{E_{(\alpha,\beta)}}}(\mathfrak{a_{1},a_{2};b_{1};}u)$ are obtained in terms of the generalized HMLF ${_{_{p}F_{q}}{E_{(\alpha,\beta)}}}(\mathfrak{a_{1},a_{2},\cdots,a}_{p};\mathfrak{b_{1},b_{2},\cdots,b}_{q};u)$ by proving the following results:
\begin{thm}
For the hypergeometric-Mittag-Leffler functions ${_{_{2}F_{1}}{E_{(\alpha,\beta)}}}(\mathfrak{a_{1},a_{2};b_{1};}u)$, the following integral representation holds true:
\begin{equation}\label{mheq31}
\int u^{\delta}{_{_{2}F_{1}}{E_{(\alpha,\beta)}}}(\mathfrak{a_{1},a_{2};b_{1};}u)du=\dfrac{u^{\delta+1}}{\delta+1}\;{_{_{3}F_{2}}{E_{(\alpha,\beta)}}}(\mathfrak{a_{1},a_{2}},\delta+1;\mathfrak{b_{1}},\delta+2;u).
\end{equation}
\begin{proof}
Multiplying equation \eqref{mheq15} by $u^{\delta}$ and integrating the resultant equation w.r.t. $u$, it follows that:
\begin{equation*}
	\int u^{\delta}{_{_{2}F_{1}}{E_{(\alpha,\beta)}}}(\mathfrak{a_{1},a_{2};b_{1};}u)du=\sum_{r=0}^{\infty}\dfrac{u^{(r+\delta+1)}{_{2}\hat{\chi}^{r}_{1}}\;\hat{d}^{r}_{(\alpha,\beta)}}{r!(\delta+r+1)}\Psi_{0}\Phi_{0}.
\end{equation*}
Applying umbral actions \eqref{mheq6} and \eqref{mheq8} in the r.h.s. of above equation and using the following identity:
\begin{equation*}
	\dfrac{1}{\delta+r+1}=\dfrac{\Gamma(\delta+1)(\delta+1)_{r}}{\Gamma(\delta+2)(\delta+2)_{r}},
\end{equation*}
we find
\begin{equation}\label{mheq32}
	\int u^{\delta}{_{_{2}F_{1}}{E_{(\alpha,\beta)}}}(\mathfrak{a_{1},a_{2};b_{1};}u)du	=\dfrac{u^{\delta+1}}{\delta+1}\sum_{r=0}^{\infty}\dfrac{(\mathfrak{a_{1}})_{r}(\mathfrak{a_{2}})_{r}(\delta+1)_{r}}{(\mathfrak{b_{1}})_{r}(\delta+2)_{r}\Gamma(\alpha r+\beta)}u^{r},
\end{equation}
which in view of definition \eqref{mheq86}, yields assertion \eqref{mheq31}.
\end{proof}
\end{thm}
\begin{rem}
Following the lines of proof of Theorem \thmnumber{5}, we obtain the analogous result for the generalized HMLF ${_{_{p}F_{q}}{E_{(\alpha,\beta)}}}(\mathfrak{a_{1}},\mathfrak{a_{2}},\cdots,\mathfrak{a}_{p};\mathfrak{b_{1}},\mathfrak{b_{2}},\cdots,\mathfrak{b}_{q};u)$.
\end{rem}
\begin{thm}
The following integral representation holds for the generalized hypergeometric-Mittag-Leffler functions ${_{_{p}F_{q}}{E_{(\alpha,\beta)}}}(\mathfrak{a_{1}},\mathfrak{a_{2}},\cdots,\mathfrak{a}_{p};\mathfrak{b_{1}},\mathfrak{b_{2}},\cdots,\mathfrak{b}_{q};u)$:
\begin{align}\label{mheq19}
\notag	\int u^{\delta}{_{_{p}F_{q}}{E_{(\alpha,\beta)}}}(\mathfrak{a_{1}},\mathfrak{a_{2}},\cdots,\mathfrak{a}_{p};\mathfrak{b_{1}},\mathfrak{b_{2}}&,\cdots,\mathfrak{b}_{q};u)du=\dfrac{u^{\delta+1}}{\delta+1}\\&\times{_{_{p+1}F_{q+1}}{E_{(\alpha,\beta)}}}(\mathfrak{a_{1},a_{2}},\cdots,\mathfrak{a_{p}},\delta+1;\mathfrak{b_{1}},\cdots,\mathfrak{b_{q}},\delta+2;u).
\end{align}
\end{thm}
\begin{rem}
Taking $\delta=0$ in equations \eqref{mheq31} and \eqref{mheq19}, we deduce the following consequences of Theorems \thmnumber{5} and \thmnumber{6}:
\end{rem}
\begin{cor}
For the hypergeometric-Mittag-Leffler functions ${_{_{2}F_{1}}{E_{(\alpha,\beta)}}}(\mathfrak{a_{1},a_{2};b_{1};}u)$ and generalized hypergeometric-Mittag-Leffler functions ${_{_{p}F_{q}}{E_{(\alpha,\beta)}}}(\mathfrak{a_{1}},\mathfrak{a_{2}},\cdots,\mathfrak{a}_{p};\mathfrak{b_{1}},\mathfrak{b_{2}},\cdots,\mathfrak{b}_{q};u)$, the following integral representations hold:
\begin{equation*}
\int{_{_{2}F_{1}}{E_{(\alpha,\beta)}}}(\mathfrak{a_{1},a_{2};b_{1};}u)du=u\;{_{_{3}F_{2}}{E_{(\alpha,\beta)}}}(\mathfrak{a_{1},a_{2}},1;\mathfrak{b_{1}},2;u),
\end{equation*}
\begin{align*}
\notag	\int{_{_{p}F_{q}}{E_{(\alpha,\beta)}}}(\mathfrak{a_{1}},\mathfrak{a_{2}},\cdots,\mathfrak{a}_{p};\mathfrak{b_{1}},\mathfrak{b_{2}},\cdots,\mathfrak{b}_{q};u)du=u\;{_{_{p+1}F_{q+1}}{E_{(\alpha,\beta)}}}(\mathfrak{a_{1},a_{2}},\cdots,\mathfrak{a_{p}},1;\mathfrak{b_{1}},\cdots,\mathfrak{b_{q}},2;u).
\end{align*}
\end{cor}
\begin{thm}
For the hypergeometric-Mittag-Leffler functions ${_{_{2}F_{1}}{E_{(\alpha,\beta)}}}(\mathfrak{a_{1},a_{2};b_{1};}u)$, the following integral representation holds:
\begin{equation}\label{mheq72}
\int_{-\infty}^{\infty}e^{-\delta u^{2}}	{_{_{2}F_{1}}{E_{(\alpha,\beta)}}}(\mathfrak{a_{1},a_{2};b_{1};}u)du=\sqrt{\dfrac{\pi}{\delta}}\;{_{_{5}F_{2}}{E_{(2\alpha,\beta)}}}\left(\dfrac{\mathfrak{a_{1}}}{2},\dfrac{\mathfrak{a_{1}}+1}{2},\dfrac{\mathfrak{a_{2}}}{2},\dfrac{\mathfrak{a_{2}}+1}{2},\dfrac{1}{2};\dfrac{\mathfrak{b_{1}}}{2},\dfrac{\mathfrak{b_{1}}+1}{2};\dfrac{4}{\delta}\right).
\end{equation}
\begin{proof}
Multiplying equation \eqref{mheq12} by $e^{-\delta u^{2}}$ and integrating w.r.t. $u$ between the limits $-\infty$ to $\infty$ by making use of integral
\begin{equation}\label{mheq91}
	\int_{-\infty}^{\infty}e^{-(ax^{2}+bx+c)}dx=\sqrt{\dfrac{\pi}{a}}\exp\left\{\dfrac{b^{2}}{4a}+c\right\},
\end{equation} 
it follows that:
\begin{equation}\label{mheq67}
	\int_{-\infty}^{\infty}e^{-\delta u^{2}}	{_{_{2}F_{1}}{E_{(\alpha,\beta)}}}(\mathfrak{a_{1},a_{2};b_{1};}u)du=\sqrt{\dfrac{\pi}{\delta}}\exp\left\{{\dfrac{\left(\hat{_{2}\chi_{1}}\;\hat{d}_{(\alpha,\beta)}\right)^{2}}{4\delta}}\right\}\Psi_{0}\Phi_{0}.
\end{equation}
Expanding the exponential in the r.h.s. of equation \eqref{mheq67} and using umbral operators \eqref{mheq6} and \eqref{mheq8}, we find
\begin{equation*}
	\int_{-\infty}^{\infty}e^{-\delta u^{2}}	{_{_{2}F_{1}}{E_{(\alpha,\beta)}}}(\mathfrak{a_{1},a_{2};b_{1};}u)du=\sqrt{\dfrac{\pi}{\delta}}\;\sum_{r=0}^{\infty}\dfrac{\Gamma(2r+1)(\mathfrak{a_{1}})_{2r}(\mathfrak{a_{2}})_{2r}}{(\mathfrak{b_{1}})_{2r}(4\delta)^{r}\Gamma(2\alpha r+\beta)\;r!},
\end{equation*}
which on making use of identities 
\begin{equation}\label{mheqxii}
	(d)_{2r}=2^{2r}\left(\dfrac{d}{2}\right)_{r}\left(\dfrac{d+1}{2}\right)_{r}
\end{equation}
and
\begin{equation}\label{mheq20}
	(2r)!=\dfrac{1}{\sqrt{\pi}}2^{2r}r!\;\Gamma(r+\dfrac{1}{2}),
\end{equation}
gives
\begin{equation*}
	\int_{-\infty}^{\infty}e^{-\delta u^{2}}	{_{_{2}F_{1}}{E_{(\alpha,\beta)}}}(\mathfrak{a_{1},a_{2};b_{1};}u)du=\sqrt{\dfrac{\pi}{\delta}}\;\sum_{r=0}^{\infty}\dfrac{\left(\dfrac{\mathfrak{a_{1}}}{2}\right)_{r}\left(\dfrac{\mathfrak{a_{1}}+1}{2}\right)_{r}\left(\dfrac{\mathfrak{a_{2}}}{2}\right)_{r}\left(\dfrac{\mathfrak{a_{2}}+1}{2}\right)_{r}\left(\dfrac{1}{2}\right)_{r}}{\left(\dfrac{\mathfrak{b_{1}}}{2}\right)_{r}\left(\dfrac{\mathfrak{b_{1}}+1}{2}\right)_{r}\Gamma(2\alpha r+\beta)}\left(\dfrac{4}{\delta}\right)^{r}.
\end{equation*}
Finally, using equation \eqref{mheq86}, assertion \eqref{mheq72} follows.
\end{proof}
\end{thm}
Following the lines of proof of Theorem \thmnumber{7}, the integral representation for the generalized HMLF ${_{_{q+1}F_{q}}{E_{(\alpha,\beta)}}}(\mathfrak{a_{1}},\mathfrak{a_{2}},\cdots,\mathfrak{a}_{q+1};\mathfrak{b_{1}},\mathfrak{b_{2}},\cdots,\mathfrak{b}_{q};u)$ is obtained by iteration method in the form of following result:
\begin{thm}
For the generalized hypergeometric-Mittag-Leffler functions ${_{_{q+1}F_{q}}{E_{(\alpha,\beta)}}}(\mathfrak{a_{1},a_{2},\cdots,a}_{q+1};\mathfrak{b_{1},b_{2},\cdots,b}_{q};u)$, the following integral holds true:
\begin{align}\label{mheq1}
\notag &\int_{-\infty}^{\infty}e^{-\delta u^{2}}	{_{_{q+1}F_{q}}{E_{(\alpha,\beta)}}}(\mathfrak{a_{1},a_{2},\cdots,a}_{q+1};\mathfrak{b_{1},b_{2},\cdots,b}_{q};u)\;du=\sqrt{\dfrac{\pi}{\delta}}\\	&\times{_{_{2(q+1)+1}F_{2q}}{E_{(2\alpha,\beta)}}}
\left(\dfrac{\mathfrak{a_{1}}}{2},\dfrac{\mathfrak{a_{1}}+1}{2},\cdots,\dfrac{\mathfrak{a}_{2(q+1)}}{2},\dfrac{\mathfrak{a}_{2(q+1)}+1}{2},\dfrac{1}{2};\dfrac{\mathfrak{b_{1}}}{2},\dfrac{\mathfrak{b_{1}}+1}{2},\cdots,\dfrac{\mathfrak{b}_{2q}}{2},\dfrac{\mathfrak{b}_{2q}+1}{2};\dfrac{4}{\delta}\right),\notag\\&\hspace{13.2cm}|u|<1.
\end{align}
\end{thm}
\begin{rem}
By extending the approach used in proving Theorems \thmnumber{7} and \thmnumber{8}, we get the following results:
\end{rem}
\begin{thm}
Following integral representation holds true for the generalized hypergeometric-Mittag-Leffler functions ${_{_{q}F_{q}}{E_{(\alpha,\beta)}}}(\mathfrak{a_{1},a_{2},\cdots,a}_{q};\mathfrak{b_{1},b_{2},\cdots,b}_{q};u)$:
\begin{align}\label{mheq2}
\notag \int_{-\infty}^{\infty}e^{-\delta u^{2}}&	{_{_{q}F_{q}}{E_{(\alpha,\beta)}}}(\mathfrak{a_{1},a_{2},\cdots,a}_{q};\mathfrak{b_{1},b_{2},\cdots,b}_{q};u)\;du=\sqrt{\dfrac{\pi}{\delta}}\\&\times{_{_{2q+1}F_{2q}}{E_{(2\alpha,\beta)}}}
\left(\dfrac{\mathfrak{a_{1}}}{2},\dfrac{\mathfrak{a_{1}}+1}{2},\cdots,\dfrac{\mathfrak{a}_{2q}}{2},\dfrac{\mathfrak{a}_{2q}+1}{2},\dfrac{1}{2};\dfrac{\mathfrak{b_{1}}}{2},\dfrac{\mathfrak{b_{1}}+1}{2},\cdots,\dfrac{\mathfrak{b}_{2q}}{2},\dfrac{\mathfrak{b}_{2q}+1}{2};\dfrac{1}{\delta}\right)
\end{align}
\end{thm}
\begin{thm}
For the generalized hypergeometric-Mittag-Leffler functions ${_{_{q}F_{q+1}}{E_{(\alpha,\beta)}}}(\mathfrak{a_{1},a_{2},\cdots,a}_{q};\mathfrak{b_{1},b_{2},\cdots,b}_{q+1};u)$, the following integral representation holds true:
\begin{align}\label{mheq3}
\notag \int_{-\infty}^{\infty}&e^{-\delta u^{2}}	{_{_{q}F_{q+1}}{E_{(\alpha,\beta)}}}(\mathfrak{a_{1},a_{2},\cdots,a}_{q+1};\mathfrak{b_{1},b_{2},\cdots,b}_{q};u)\;du=\sqrt{\dfrac{\pi}{\delta}}\\	&\times{_{_{2q+1}F_{2(q+1)}}{E_{(2\alpha,\beta)}}}
\left(\dfrac{\mathfrak{a_{1}}}{2},\dfrac{\mathfrak{a_{1}}+1}{2},\cdots,\dfrac{\mathfrak{a}_{2q}}{2},\dfrac{\mathfrak{a}_{2q}+1}{2},\dfrac{1}{2};\dfrac{\mathfrak{b_{1}}}{2},\dfrac{\mathfrak{b_{1}}+1}{2},\cdots,\dfrac{\mathfrak{b}_{2(q+1)}}{2},\dfrac{\mathfrak{b}_{2(q+1)}+1}{2};\dfrac{1}{4\delta}\right).
\end{align}
\end{thm}
Finally, we derive the integral representations involving the product of the HMLF ${_{_{2}F_{1}}{E_{(\alpha,\beta)}}}(\mathfrak{a_{1},a_{2};b_{1};}u)$ with $sine$ and $cosine$ functions.
\begin{thm}
Following integral representations involving the product of hypergeometric-Mittag-Leffler functions ${_{_{2}F_{1}}{E_{(\alpha,\beta)}}}(\mathfrak{a_{1},a_{2};b_{1};}u)$ with $sine$ and $cosine$ functions hold true:
\begin{equation}\label{mheq77}
\int_{0}^{\infty}{_{_{2}F_{1}}{E_{(\alpha,\beta)}}}(\mathfrak{a_{1},a_{2};b_{1};}-u)\sin(u)du={_{_{6}F_{2}}{E_{(2\alpha,\beta)}}}\left(\dfrac{\mathfrak{a_{1}}}{2},\dfrac{\mathfrak{a_{1}}+1}{2},\dfrac{\mathfrak{a_{2}}}{2},\dfrac{\mathfrak{a_{2}}+1}{2},\dfrac{1}{2},1;\dfrac{\mathfrak{b_{1}}}{2},\dfrac{\mathfrak{b_{1}}+1}{2};-16\right)
\end{equation}
and
\begin{align}\label{mheq76}
\notag
\int_{0}^{\infty}{_{_{2}F_{1}}{E_{(\alpha,\beta)}}}(\mathfrak{a_{1},a_{2};}&\mathfrak{b_{1};}-u)\cos(u)du=\dfrac{\mathfrak{a_{1}a_{2}}}{\mathfrak{b_{1}}}\\ &\times{_{_{6}F_{2}}{E_{(2\alpha,\alpha+\beta)}}}\left(\dfrac{\mathfrak{a_{1}}+1}{2},\dfrac{\mathfrak{a_{1}}+2}{2},\dfrac{\mathfrak{a_{2}}+1}{2},\dfrac{\mathfrak{a_{2}}+2}{2},\dfrac{3}{2},1;\dfrac{\mathfrak{b_{1}}+1}{2},\dfrac{\mathfrak{b_{1}}+2}{2};-16\right),
\end{align}
respectively.
\begin{proof}
Replacing $u$ by $-u$ in equation \eqref{mheq12} then multiplying the result by $e^{-iu}$ and integrating both sides w.r.t. $u$ between the limits $0$ and $\infty$, we have
\begin{equation}\label{mheq79}
	\int_{0}^{\infty}{_{_{2}F_{1}}{E_{(\alpha,\beta)}}}(\mathfrak{a_{1},a_{2};b_{1};}-u)e^{-iu}du=\dfrac{-i}{1-i\;{_{2}\hat{\chi}_{1}}\;\hat{d}_{(\alpha,\beta)}}\Psi_{0}\Phi_{0}.
\end{equation}
Performing certain rearrangements in equation \eqref{mheq79} and applying operator equations \eqref{mheq6} and \eqref{mheq8}, it follows that:
\begin{equation*}
	\int_{0}^{\infty}{_{_{2}F_{1}}{E_{(\alpha,\beta)}}}(\mathfrak{a_{1},a_{2};b_{1};}-u)e^{-iu}du=-i\;\sum_{r=0}^{\infty}i^{r}\hat{_{2}\chi^{r}_{1}}\;\hat{d}^{r}_{(\alpha,\beta)}\Psi_{0}\Phi_{0}=-i\;\sum_{r=0}^{\infty}\dfrac{(\mathfrak{a}_{1})_{r}(\mathfrak{a}_{2})_{r}(1)_{r}}{(\mathfrak{b}_{1})_{r}\Gamma(\alpha r+\beta)}i^{r},
\end{equation*}
which in view of equations \eqref{mheqv} and \eqref{mheq86} takes the form
\begin{equation}\label{mheq80}
	\int_{0}^{\infty}{_{_{2}F_{1}}{E_{(\alpha,\beta)}}}(\mathfrak{a_{1},a_{2};b_{1};}-u)e^{-iu}du=-i\;e^{i\;{_{3}\hat{\chi}_{1}}\;\hat{d}_{(\alpha,\beta)}}\Psi_{0}\Phi_{0},
\end{equation}
Consequently, we have
\begin{equation}\label{mheq78}
	\int_{0}^{\infty}{_{_{2}F_{1}}{E_{(\alpha,\beta)}}}(\mathfrak{a_{1},a_{2};b_{1};}-u)\left(\cos u-i\sin u\right)du=\left\{\sin\left({_{3}\hat{\chi}_{1}}\;\hat{d}_{(\alpha,\beta)}\right)-i\cos\left({_{3}\hat{\chi}_{1}}\;\hat{d}_{(\alpha,\beta)}\right)\right\}\Psi_{0}\Phi_{0}.
\end{equation}
Equating imaginary parts in equation \eqref{mheq78}, we obtain
\begin{equation*}
	\int_{0}^{\infty}{_{_{2}F_{1}}{E_{(\alpha,\beta)}}}(\mathfrak{a_{1},a_{2};b_{1};}-u)\sin(u)du=\cos\left({_{3}\hat{\chi}_{1}}\;\hat{d}_{(\alpha,\beta)}\right)\Psi_{0}\Phi_{0},
\end{equation*}
which on using series expansion of $cosine~function$ and applying operator actions \eqref{mheqii} and \eqref{mheq8}, becomes
\begin{equation*}
	\int_{0}^{\infty}{_{_{2}F_{1}}{E_{(\alpha,\beta)}}}(\mathfrak{a_{1},a_{2};b_{1};}-u)\sin(u)du=\sum_{r=0}^{\infty}\dfrac{(-1)^{r}(\mathfrak{a_{1}})_{2r}(\mathfrak{a_{2}})_{2r}(1)_{2r}\Gamma(2r+1)}{(\mathfrak{b_{1}})_{2r}(2r)!\;\Gamma(2\alpha r+\beta)}.
\end{equation*}
Use of identity \eqref{mheqxii} in the above equation yields assertion \eqref{mheq77}.\\

Further, equating the real parts in equation \eqref{mheq78}, we find
\begin{equation*}
	\int_{0}^{\infty}{_{_{2}F_{1}}{E_{(\alpha,\beta)}}}(\mathfrak{a_{1},a_{2};b_{1};}-u)\cos(u)du=\sin\left({_{3}\hat{\chi}_{1}}\;\hat{d}_{(\alpha,\beta)}\right)\Psi_{0}\Phi_{0}.
\end{equation*}
Expanding the $sine~function$ and applying umbral actions \eqref{mheqii} and \eqref{mheq8}, we have
\begin{equation*}
	\int_{0}^{\infty}{_{_{2}F_{1}}{E_{(\alpha,\beta)}}}(\mathfrak{a_{1},a_{2};b_{1};}-u)\cos(u)du=\sum_{r=0}^{\infty}\dfrac{(-1)^{r}(\mathfrak{a_{1}})_{2r+1}(\mathfrak{a_{2}})_{2r+1}(1)_{2r+1}\Gamma(2r+2)}{(\mathfrak{b_{1}})_{2r+1}(2r+1)!\;\Gamma(2\alpha r+(\alpha+\beta))},
\end{equation*}
which on making use of identities \eqref{mheqxi} and \eqref{mheqxii} yields assertion \eqref{mheq76}.
\end{proof}
\end{thm}
\noindent In the next section, the Laplace and Sumudu transforms for the hypergeometric-Mittag-Leffler functions are established. The graphical representations and distribution of zeros of the HMLF ${_{_{2}F_{1}}{E_{(\alpha,\beta)}}}(\mathfrak{a_{1},a_{2};b_{1};}u)$ are also considered.
\section{Conclusions}
We have studied many facets of the umbral formalism and its application to the HMLF in previous sections. The distinctive characteristics of the stated class of functions are explored and analyzed with special emphasis on the relevant series representation, differential and integral formulae.\\ 

Integral transformations and special functions play an important role in mathematics and engineering. Here, we extend the formalism to explore the Laplace and Sumudu transformation. The Laplace and Sumudu transformations have given a new spirit in solving several problems confronted in engineering applications.\\

We recall that the Laplace transform \cite{MRS} of the function $f(t)$ is defined as:
\begin{equation}\label{mheqxvi}
F(s)=\mathcal{L}{\{f(t):s\}}=\int_{0}^{\infty}e^{-st}f(t)dt,\quad\quad\quad t>0,
\end{equation}
provided the limit exists.\\

Also, we recall that the Sumudu transform \cite{GKW} is defined as:
\begin{equation}\label{mheqxv}
G(u)=\mathcal{S}{\{f(t): u\}}=\int_{0}^{\infty}e^{-t}f(ut)dt.
\end{equation}
Replacing $u$ by $-t$ in equation \eqref{mheq12} then multiplying the result by $e^{-st}$ and integrating both sides w.r.t. $t$ between the limits $0$ and $\infty$, we have
\begin{equation*}
\int_{0}^{\infty}e^{-st}{_{_{2}F_{1}}{E_{(\alpha,\beta)}}}(\mathfrak{a_{1},a_{2};b_{1};}-t)dt=\int_{0}^{\infty}e^{-st-t\;{_{2}\hat{\chi}_{1}}\;\hat{d}_{(\alpha,\beta)}}\Psi_{0}\Phi_{0}dt,
\end{equation*}
which on performing integration in the r.h.s. takes the form
\begin{equation}\label{mheq36}
\int_{0}^{\infty}e^{-st}{_{_{2}F_{1}}{E_{(\alpha,\beta)}}}(\mathfrak{a_{1},a_{2};b_{1};}-t)dt=\dfrac{1}{s+{_{2}\hat{\chi}_{1}}\;\hat{d}_{(\alpha,\beta)}}\Psi_{0}\Phi_{0}.
\end{equation}
Performing certain rearrangements in equation \eqref{mheq36} and applying operators \eqref{mheq6} and \eqref{mheq8} in the resultant, it follows that:
\begin{equation*}
\int_{0}^{\infty}e^{-st}{_{_{2}F_{1}}{E_{(\alpha,\beta)}}}(\mathfrak{a_{1},a_{2};b_{1};}-t)dt	=\dfrac{1}{s}\sum_{r=0}^{\infty}\dfrac{(\mathfrak{a}_{1})_{r}(\mathfrak{a}_{2})_{r}(1)_{r}}{(\mathfrak{b}_{1})_{r}\Gamma(\alpha r+\beta)}\left(\dfrac{-1}{s}\right)^{r}.
\end{equation*}
In view of equations \eqref{mheqv}, \eqref{mheq86} and definition \eqref{mheqxvi}, we get the following Laplace transform for HMLF ${_{_{2}F_{1}}{E_{(\alpha,\beta)}}}(\mathfrak{a_{1},a_{2};b_{1};}u)$:
\begin{equation}\label{mheq44}
\mathcal{L}\left\{{_{_{2}F_{1}}{E_{(\alpha,\beta)}}}(\mathfrak{a_{1},a_{2};b_{1};}-t):s\right\}=\dfrac{1}{s}\;{_{_{3}F_{1}}{E_{(\alpha,\beta)}}}(\mathfrak{a_{1},a_{2},1;b_{1};}-1/s).
\end{equation}

Again, replacing $u$ by $-ut$ in equation \eqref{mheq12} then multiplying by $e^{-t}$ and integrating w.r.t. $t$ between the limits $0$ and $\infty$, we have
\begin{equation*}
\int_{0}^{\infty}e^{-t}{_{_{2}F_{1}}{E_{(\alpha,\beta)}}}(\mathfrak{a_{1},a_{2};b_{1};}-ut)dt=\dfrac{1}{1+u\;{_{2}\hat{\chi}_{1}}\;\hat{d}_{(\alpha,\beta)}}\Psi_{0}\Phi_{0}.
\end{equation*}
Following the lines of proof of transform \eqref{mheq44} and using definition \eqref{mheqxv} in the resultant, yields the following Sumudu transformation for HMLF ${_{_{2}F_{1}}{E_{(\alpha,\beta)}}}(\mathfrak{a_{1},a_{2};b_{1};}u)$:
\begin{equation}\label{mheq46}
\mathcal{S}\left\{{_{_{2}F_{1}}{E_{(\alpha,\beta)}}}(\mathfrak{a_{1},a_{2};b_{1};}-t):u\right\}={_{_{3}F_{1}}{E_{(\alpha,\beta)}}}(\mathfrak{a_{1},a_{2},1;b_{1};}-u).
\end{equation}
Graphical representations are valuable assets across numerous domains such as science, engineering, technology, finance, and other fields. They offer swift and intuitive insights into the behaviour of functions through graphs and plots. Consequently, this increases the effectiveness of the research process.\\

We aim to explore the geometry of HMLF ${_{_{2}F_{1}}{E_{(\alpha,\beta)}}}(\mathfrak{a_{1},a_{2};b_{1};}u)$ using the software \textquotedblleft Mathematica\textquotedblright. Using series expansion \eqref{mheq14}, we draw the graphs of HMLF ${_{_{2}F_{1}}{E_{(\alpha,\beta)}}}(\mathfrak{a_{1},a_{2};b_{1};}u)$.\\

Figure 1, represents the graph of HMLF ${_{_{2}F_{1}}{E_{(2,3)}}}(3,1;3;u)$.\\
Figure 2, represents the graph of HMLF ${_{_{2}F_{1}}{E_{(2,1)}}}(1,1;1;u)$.\\
Figure 3, represents the graph of HMLF ${_{_{2}F_{1}}{E_{(1,2)}}}(3,5;1;u)$.\\
Figure 4, represents the graph of HMLF ${_{_{2}F_{1}}{E_{(1,4)}}}(1,2;3;u)$.\\
\begin{figure}[h!]
\begin{minipage}[b]{.44\textwidth}
\includegraphics[width=0.9\linewidth]{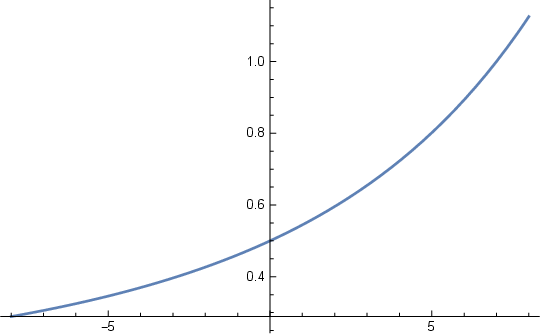}
\caption{Graph of ${_{_{2}F_{1}}{E_{(2,3)}}}(3,1;3;u)$}
\end{minipage}
\hspace{1cm}
\begin{minipage}[b]{.44\textwidth}
\includegraphics[width=0.9\linewidth]{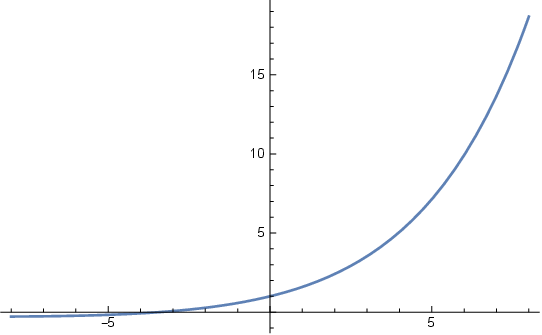}
\caption{Graph of ${_{_{2}F_{1}}{E_{(2,1)}}}(1,1;1;u)$}
\end{minipage}
\end{figure}
\begin{figure}[h!]
\begin{minipage}[b]{.44\textwidth}
\includegraphics[width=0.9\linewidth]{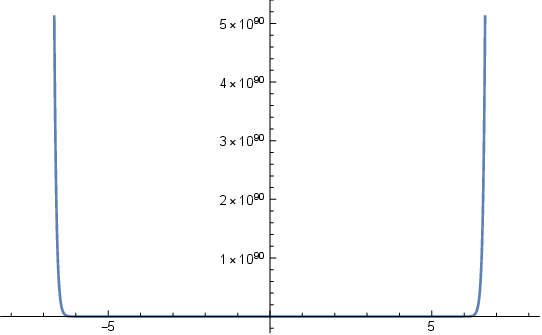}
\caption{Graph of ${_{_{2}F_{1}}{E_{(1,2)}}}(3,5;1;u)$}
\end{minipage}
\hspace{1cm}
\begin{minipage}[b]{.44\textwidth}
\includegraphics[width=0.9\linewidth]{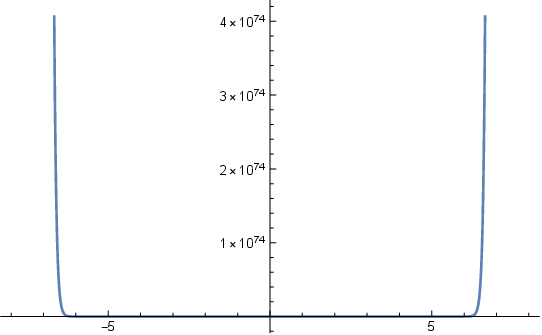}
\caption{Graph of ${_{_{2}F_{1}}{E_{(1,4)}}}(1,2;3;u)$}
\end{minipage}
\end{figure}

\noindent To elucidate the distribution of pattern of zeros of the HMLF ${_{_{2}F_{1}}{E_{(\alpha,\beta)}}}(\mathfrak{a_{1},a_{2};b_{1};}u)$, we visually represent these zeros for specific values of the parameters.
\begin{figure}[h!]
\begin{minipage}[b]{.44\textwidth}
\includegraphics[width=0.9\linewidth]{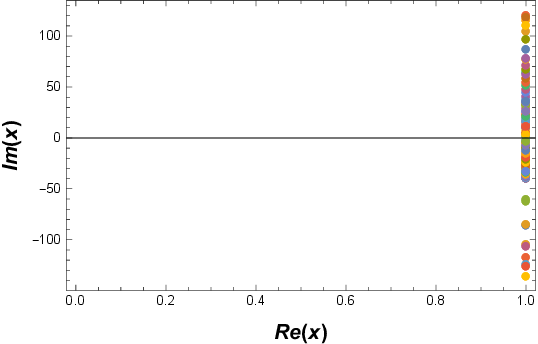}
\caption{Zeros of ${_{_{2}F_{1}}{E_{(1,3)}}}(1.5,2;1.5;u)$}
\end{minipage}
\hspace{1cm}
\begin{minipage}[b]{.44\textwidth}
\includegraphics[width=0.9\linewidth]{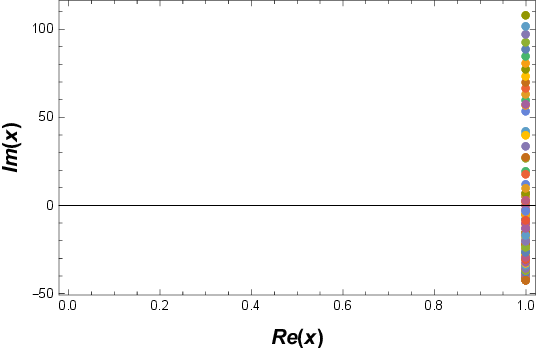}
\caption{Graph of ${_{_{2}F_{1}}{E_{(1,4)}}}(1,1;1;u)$}
\end{minipage}
\end{figure}
These illustrations show different pattern of distribution of zeros of HMLF ${_{_{2}F_{1}}{E_{(\alpha,\beta)}}}(\mathfrak{a_{1},a_{2};b_{1};}u)$.\\

Figure 5, represents the zeros of HMLF ${_{_{2}F_{1}}{E_{(1,3)}}}(1.5,2;1.5;u)$.\\
Figure 6, represents the zeros of HMLF ${_{_{2}F_{1}}{E_{(1,4)}}}(1,1;1;u)$.\\ 

\noindent We observe that until now, our discussion has been focused on the single-variable HMLF, which has significant relevance in various contexts. However, prior to concluding this article, it is important to note that multi-variable extension may be considered, especially when exploring practical applications. In an upcoming study, we plan to utilize and expand upon the formalism, we have developed so far, aiming to establish the umbral restyling of multi-variable HMLF and other hybrid forms. These extensions will allow us to delve deeper into the complexities and broader implications of these functions in diverse mathematical and scientific domains. 
\section*{Appendix}
For suitable selection of parameters in the generalized hypergeometric-Mittag-Leffler functions ${_{_{p}F_{q}}{E_{(\alpha,\beta)}}}(\mathfrak{a_{1},a_{2},\cdots,a}_{p};\mathfrak{b_{1},b_{2},\cdots,b}_{q};u)$, we get certain important connections of HMLF with other functions.\\
{\bf 1.} For $p=4$, $q=3$; taking $\alpha=1$ and $\beta=m+1$ in equation \eqref{mheq86}, the following relation is\newline\hspace*{0.4cm} obtained:
\begin{align}
\notag
{_{_{4}F_{3}}{E_{(1,m+1)}}}\left(\dfrac{\mathfrak{a_{1}}+m}{2},\dfrac{\mathfrak{a_{1}}+m+1}{2},\dfrac{\mathfrak{a_{2}}+m}{2},\dfrac{\mathfrak{a_{2}}+m+1}{2};1,\dfrac{\mathfrak{b_{1}}+m}{2},\dfrac{\mathfrak{b_{1}}+m+1}{2};-u^{2}\right)&\\\tag{A.1}=\left(\dfrac{2}{u}\right)^{m}\dfrac{\mathfrak{(b_{1})_{m}}}{\mathfrak{(a_{1})_{m}(a_{2})_{m}}}{_{_{2}F_{1}}{\mathcal{J}_{m}}}(\mathfrak{a_{1},a_{2};b_{1};}u),
\end{align}
where, hypergeometric-Bessel functions (HBF) ${_{_{2}F_{1}}{\mathcal{J}_{m}}}(\mathfrak{a_{1},a_{2};b_{1};}u)$ are recently introduced in \cite{SM}:
\begin{equation*}
{_{_{2}F_{1}}{\mathcal{J}_{m}}}(\mathfrak{a_{1},a_{2};b_{1};}u)=\sum_{r=0}^{\infty}\dfrac{(-1)^{r}(\mathfrak{a_{1}})_{m+2r}(\mathfrak{a_{2}})_{m+2r}}{r!\;\Gamma(m+r+1)(\mathfrak{b_{1}})_{m+2r}}\left(\dfrac{u}{2}\right)^{m+2r}.
\end{equation*}
Using relation (A.1) in integral representation \eqref{mheq1}, we get the following integral for HBF ${_{_{2}F_{1}}{\mathcal{J}_{m}}}(\mathfrak{a_{1},a_{2};b_{1};}u)$:
\begin{align*}
\notag &\int_{-\infty}^{\infty}\dfrac{e^{-\delta u^{2}}}{u^{m/2}}\;	{_{_{2}F_{1}}{\mathcal{J}_{m}}}(\mathfrak{a_{1},a_{2};b_{1};}\sqrt{u})\;du=\sqrt{\dfrac{\pi}{\delta}}\dfrac{\mathfrak{(a_{1})_{m}(a_{2})_{m}}}{2^{m}\mathfrak{(b_{1})_{m}}}\\\tag{A.2}&\times{_{_{8}F_{5}}{E_{(2,m+1)}}}
\left(\dfrac{\mathfrak{a_{1}}+m}{4},\dfrac{\mathfrak{a_{1}}+m+1}{4},\dfrac{\mathfrak{a_{1}}+m+2}{4},\dfrac{\mathfrak{a_{1}}+m+3}{4},\cdots
;1,\dfrac{\mathfrak{b_{1}}+m}{4},\dfrac{\mathfrak{b_{1}}+m+1}{4}\cdots;\dfrac{4}{\delta}\right).
\end{align*}
{\bf 2.} For $p=2$, $q=2$; taking $\alpha=1$ and $\beta=m+1$ in equation \eqref{mheq86}, the following connection is\newline\hspace*{0.4cm} obtained:
\begin{equation}\tag{A.3}
{_{_{2}F_{2}}{E_{(1,m+1)}}}(\mathfrak{a_{1},a_{2};b_{1}},1;-u)={_{_{2}F_{1}}{\mathcal{C}_{m}}}(\mathfrak{a_{1},a_{2};b_{1};}u),
\end{equation}
where, the hypergeometric-Tricomi functions (HTF) ${_{_{2}F_{1}}{\mathcal{C}_{m}}}(\mathfrak{a_{1},a_{2};b_{1};}u)$ are given as \cite{SM}:
\begin{equation*}
{_{_{2}F_{1}}{\mathcal{C}_{m}}}(\mathfrak{a_{1},a_{2};b_{1};}u)=\sum_{r=0}^{\infty}\dfrac{(-u)^{r}(\mathfrak{a}_{1})_{r}(\mathfrak{a}_{2})_{r}}{\Gamma(m+r+1)\;r!\;(\mathfrak{b}_{1})_{r}}.
\end{equation*}
Making use of connection (A.3) in integrals \eqref{mheq19} and \eqref{mheq2}, the following integrals for HTF ${_{_{2}F_{1}}{\mathcal{C}_{m}}}(\mathfrak{a_{1},a_{2};b_{1};}u)$ are obtained:
\begin{align}\tag{A.4}
\int u^{\delta}{_{_{2}F_{1}}{\mathcal{C}_{m}}}(\mathfrak{a_{1},a_{2};b_{1};}-u)du=\dfrac{u^{\delta+1}}{\delta+1}{_{_{3}F_{3}}{E_{(1,m+1)}}}(\mathfrak{a_{1},a_{2}},\delta+1;\mathfrak{b_{1}},1,\delta+2;u),
\end{align}
\begin{align*}
\notag	\int_{-\infty}^{\infty}e^{-\delta u^{2}}	{_{_{2}F_{1}}{\mathcal{C}_{m}}}(\mathfrak{a_{1},a_{2};b_{1};}-u)\;&du	=\sqrt{\dfrac{\pi}{\delta}}\\\tag{A.5}&\times{_{_{4}F_{3}}{E_{(2,m+1)}}}
\left(\dfrac{\mathfrak{a_{1}}}{2},\dfrac{\mathfrak{a_{1}}+1}{2},\dfrac{\mathfrak{a}_{2}}{2},\dfrac{\mathfrak{a}_{2}+1}{2};\dfrac{\mathfrak{b_{1}}}{2},\dfrac{\mathfrak{b_{1}}+1}{2},1;\dfrac{1}{\delta}\right).
\end{align*}
{\bf 3.} For $p=2$, $q=1$; taking $b_{1}=1$ in equation \eqref{mheq86}, we find the following relation between the\newline\hspace*{0.4cm} HMLF ${_{_{2}F_{1}}{E_{(\alpha,\beta)}}}(\mathfrak{a_{1},a_{2};}1;u)$ and Fox-Wright function ${_{2}\psi_{1}}$ \cite{AND}:
\begin{equation}\tag{A.6}
{_{_{2}F_{1}}{E_{(\alpha,\beta)}}}(\mathfrak{a_{1},a_{2};}1;u)=\dfrac{1}{\Gamma(\mathfrak{a_{1}})\Gamma(\mathfrak{a_{2}})}\;{_{2}\psi_{1}}\left[\begin{matrix}
(\mathfrak{a_{1}},1) & & (\mathfrak{a_{2}},1)\\ 
& (\alpha,\beta) &
\end{matrix};u\right].
\end{equation}
Making use of relation (A.6) in equations \eqref{mheq31} and \eqref{mheq72}, we get following results for Fox-Wright function ${_{2}\psi_{1}}$:
\begin{equation}\tag{A.7}
\int u^{\delta}{_{2}\psi_{1}}\left[\begin{matrix}
(\mathfrak{a_{1}},1) & & (\mathfrak{a_{2}},1)\\ 
& (\alpha,\beta) &
\end{matrix};u\right]du=\dfrac{\Gamma(\mathfrak{a_{1}})\Gamma(\mathfrak{a_{2}})u^{\delta+1}}{\delta+1}\;{_{_{3}F_{2}}{E_{(\alpha,\beta)}}}(\mathfrak{a_{1},a_{2}},\delta+1;1,\delta+2;u),
\end{equation}
\begin{align*}
\dfrac{1}{\Gamma(\mathfrak{a_{1}})\Gamma(\mathfrak{a_{2}})}\int_{-\infty}^{\infty}e^{-\delta u^{2}}	{_{2}\psi_{1}}\left[\begin{matrix}
(\mathfrak{a_{1}},1) & & (\mathfrak{a_{2}},1)\\ 
& (\alpha,\beta) &
\end{matrix};u\right]&du=\sqrt{\dfrac{\pi}{\delta}}\\\tag{A.8}&\times{_{_{4}F_{1}}{E_{(2\alpha,\beta)}}}\left(\dfrac{\mathfrak{a_{1}}}{2},\dfrac{\mathfrak{a_{1}}+1}{2},\dfrac{\mathfrak{a_{2}}}{2},\dfrac{\mathfrak{a_{2}}+1}{2};1;\dfrac{4}{\delta}\right).
\end{align*}
{\bf 4.} For $p=q=1$; taking $b_{1}=1$ in eqaution \eqref{mheq86}, we find the following relation between\newline\hspace*{0.4cm} confluent HMLF ${_{_{1}F_{1}}{E_{(\alpha,\beta)}}}(\mathfrak{a_{1}};1;u)$ and Prabhakar function $E_{(\alpha,\beta)}^{\;\mathfrak{a_{1}}}(u)$ \cite{P}:
\begin{equation}\tag{A.9}
{_{_{1}F_{1}}{E_{(\alpha,\beta)}}}(\mathfrak{a_{1}};1;u)=E_{(\alpha,\beta)}^{\;\mathfrak{a_{1}}}(u).
\end{equation}
Using connection (A.9) in integrals \eqref{mheq19} and \eqref{mheq2}, the following integrals are obtained for Prabhakar function $E_{(\alpha,\beta)}^{\;\mathfrak{a_{1}}}(u)$:
\begin{align}\tag{A.10}
\int u^{\delta}\;E_{(\alpha,\beta)}^{\;\mathfrak{a_{1}}}(u)du=\dfrac{u^{\delta+1}}{\delta+1}	{_{_{2}F_{2}}{E_{(\alpha,\beta)}}}(\mathfrak{a_{1}},\delta+1;1,\delta+2;u),
\end{align}
\begin{equation}\tag{A.11}
\int_{-\infty}^{\infty}e^{-\delta u^{2}}E_{(\alpha,\beta)}^{\;\mathfrak{a_{1}}}(u)du=\sqrt{\dfrac{\pi}{\delta}}{_{_{2}F_{1}}{E_{(2\alpha,\beta)}}}\left(\dfrac{\mathfrak{a_{1}}}{2},\dfrac{\mathfrak{a_{1}}+1}{2};1;\dfrac{1}{\delta}\right).
\end{equation}
{\bf 5.}	For $p=2$, $q=1$; taking $\alpha=\beta=1$ in equation \eqref{mheq86}, we get
\begin{equation}\tag{A.12}
{_{_{2}F_{1}}{E_{(1,1)}}}(\mathfrak{a_{1},a_{2};b_{1};}u)={_{2}F_{1}}(\mathfrak{a_{1},a_{2};b_{1};}u).
\end{equation} 
{\bf 6.} For $p=0$, $q=1$; taking $\alpha=b_{1}=1$ and $\beta=3/2$ in equation \eqref{mheq86}, we find
\begin{equation}\tag{A.13}
\dfrac{u\sqrt{\pi}}{2}\;{_{_{0}F_{1}}{E_{(1,3/2)}}}\left(-;1;-\dfrac{u^{2}}{4}\right)=\sin u.
\end{equation}
{\bf 7.} For $p=0$, $q=1$; taking $\alpha=b_{1}=1$ and $\beta=1/2$ in equation \eqref{mheq86}, we get
\begin{equation}\tag{A.14}
\sqrt{\pi}\;{_{_{0}F_{1}}{E_{(1,1/2)}}}\left(-;1;-\dfrac{u^{2}}{4}\right)=	\cos u.
\end{equation}
In Section $3$, the integral representations for HMLF ${_{2}F_{1}}(\mathfrak{a_{1},a_{2};b_{1};}u)$ and generalized HMLF ${_{_{p}F_{q}}{E_{(\alpha,\beta)}}}(\mathfrak{a_{1}},\mathfrak{a_{2}},\cdots,\mathfrak{a}_{p};\mathfrak{b_{1}},\mathfrak{b_{2}},\cdots,\mathfrak{b}_{q};u)$, ${_{_{q}F_{q}}{E_{(\alpha,\beta)}}}(\mathfrak{a_{1}},\mathfrak{a_{2}},\cdots,\mathfrak{a}_{q};\mathfrak{b_{1}},\mathfrak{b_{2}},\cdots,\mathfrak{b}_{q};u)$, ${_{_{q+1}F_{q}}{E_{(\alpha,\beta)}}}(\mathfrak{a_{1}},\mathfrak{a_{2}},\cdots,\mathfrak{a}_{q+1};\mathfrak{b_{1}},\mathfrak{b_{2}},\cdots,\mathfrak{b}_{q};u)$ and ${_{_{q}F_{q+1}}{E_{(\alpha,\beta)}}}(\mathfrak{a_{1}},\mathfrak{a_{2}},\cdots,\mathfrak{a}_{q};\mathfrak{b_{1}},\mathfrak{b_{2}},\cdots,\mathfrak{b}_{q+1};u)$ are established. Here, in Appendix, we obtained certain integral representations by considering special cases of these functions. However, it is not necessary that the integrals for all  these particular cases can be obtained from the results established in Section $3$. As an illustration, we derive the following integral representation for cosine function in terms of the HMLF ${_{_{0}F_{1}}{E_{(\alpha,\beta)}}}(-;\mathfrak{b_{1}};u)$:
\begin{equation}\tag{A.15}
\int_{-\infty}^{\infty}e^{-\delta u^{2}}\cos\sqrt{u}\;du=\dfrac{\pi}{\sqrt{\delta}}\;{_{_{0}F_{1}}{E_{(2,1/2)}}}(-;1;1/(64\delta)).
\end{equation}
In order to illustrate the process, we use definition \eqref{mheqv} in connection (A.15), so that we have
\begin{equation}\tag{A.16}
e^{-\delta u^{2}}\cos(\sqrt{u})=\sqrt{\pi}e^{-\delta u^{2}}{_{_{0}F_{1}}{E_{(1,1/2)}}}(-;1;-u/4)=\sqrt{\pi}e^{-\delta u^{2}}e^{(-u/4)\;\hat{_{0}\chi_{1}}\;\hat{d}_{(1,1/2)}}\Psi_{0}\Phi_{0},
\end{equation}
where
\begin{equation}\tag{A.17}
\hat{_{0}\chi^{r}_{1}}=\dfrac{1}{(1)_{r}}.
\end{equation}
Integrating (A.16) w.r.t. $u$ between the limits $-\infty$ to $\infty$ by making use of integral \eqref{mheq91},
it follows that:
\begin{equation}\tag{A.18}
\int_{-\infty}^{\infty}e^{-\delta u^{2}}\cos\sqrt{u}\;du={\dfrac{\pi}{\sqrt{\delta}}}\exp\left\{{\dfrac{\left(\hat{_{0}\chi_{1}}\;\hat{d}_{(1,1/2)}\right)^{2}}{(4)^{2}(4\delta)}}\right\}\Psi_{0}\Phi_{0}.
\end{equation}
Expanding the exponential in the r.h.s. of equation (A.18) and using umbral operators \eqref{mheq8} and (A.17), we find
\begin{equation*}
\int_{-\infty}^{\infty}e^{-\delta u^{2}}\cos\sqrt{u}\;du=\dfrac{\pi}{\sqrt{\delta}}\sum_{r=0}^{\infty}\dfrac{\Gamma(2r+1)}{(1)_{2r}\Gamma(2r+1/2)\;r!}\left(\dfrac{1}{64\delta}\right)^{r},
\end{equation*}
which on making use of identities \eqref{mheqxii} and \eqref{mheq20}, gives
\begin{equation*}
\int_{-\infty}^{\infty}e^{-\delta u^{2}}\cos\sqrt{u}\;du=\dfrac{\pi}{\sqrt{\delta}}\sum_{r=0}^{\infty}\dfrac{1}{(1)_{r}\Gamma(2r+1/2)}\left(\dfrac{1}{64\delta}\right)^{r}.
\end{equation*}
Finally, use of equation \eqref{mheq86} for $p=0$ and $q=1$, yields integral (A.15).\\

{\bf Funding}\\

This work has been done under Junior Research Fellowship (Ref No. 231610072319, dated:29/07/2023) awarded to the second author by University Grand Commission, Government of India, New Delhi.\\

{\bf Declaration of competing interest}\\

None.

\end{document}